\documentclass[12pt]{article}

\usepackage{graphicx}

\def\sm{\setminus}
\def\ra{\rightarrow}
\def\sm{\setminus}

\def\e{\epsilon}

\def\cara{Carath\'{e}odory}
\def\Aut{\hbox{Aut}}

 \def\HollowBox #1#2{{\dimen0=#1 \advance\dimen0 by -#2       
       \dimen1=#1 \advance\dimen1 by #2                       
        \vrule height #1 depth #2 width #2                    
        \vrule height 0pt depth #2 width #1                   
        \llap{\vrule height #1 depth -\dimen0 width \dimen1}%
       \hskip -#2                                             
       \vrule height #1 depth #2 width #2}}                   
 \def\BoxOpTwo{\mathord{\HollowBox{6pt}{.4pt}}\;}             

\def\endpf{\hfill $\BoxOpTwo$}

\font\teneufm=eufm10
\font\seveneufm=eufm7
\font\fiveeufm=eufm5
\newfam\eufmfam
\textfont\eufmfam=\teneufm
\scriptfont\eufmfam=\seveneufm
\scriptscriptfont\eufmfam=\fiveeufm

\newfam\msbfam
\font\tenmsb=msbm10  scaled \magstep1 \textfont\msbfam=\tenmsb
\font\sevenmsb=msbm7 scaled \magstep1 \scriptfont\msbfam=\sevenmsb
\font\fivemsb=msbm5  scaled \magstep1 \scriptscriptfont\msbfam=\fivemsb
\def\Bbb{\fam\msbfam \tenmsb}

\def\CC{{\Bbb C}}

\newtheorem{theorem}{Theorem}
\newtheorem{corollary}[theorem]{Corollary}
\newtheorem{proposition}[theorem]{Proposition}
\newtheorem{lemma}[theorem]{Lemma}

\newtheorem{definition}{Definition}
\newtheorem{example}[definition]{EXAMPLE}

\makeindex

\begin{document}

\begin{center}
\huge \bf The Carath\'{e}odory and Kobayashi Metrics 
\medskip \\
\huge \bf and Applications in Complex Analysis
\end{center}

\begin{center}
Steven G. Krantz
\end{center}

\begin{quote}
{\bf Abstract:} \sl  The Carath\'{e}odory and Kobayashi metrics have proved
to be important tools in the function theory of
several complex variables.  But they are less familiar
in the context of one complex variable.  Our purpose here is to gather
in one place the basic ideas about these important invariant
metrics for domains in the plane and to provide some illuminating examples and
applications.
\end{quote}

\setcounter{section}{-1}

\section{Prefatory Thoughts}

In the late nineteenth century, Henri Poincar\'{e} (1854--1912)
introduced the profoundly original idea of equipping the 
unit disc $D$ in the complex plane with a metric that is invariant
under conformal self-maps of $D$.  One may recall (see [GRK]) that
the conformal maps of the disc are generated by the rotations
$$
\rho_\theta : \zeta \longmapsto e^{i\theta} \zeta
$$
for $0 \leq \theta < 2\pi$
and the M\"{o}bius transformations
$$
\varphi_a : \zeta \longmapsto \frac{\zeta - a}{1 - \overline{a} \zeta}
$$
for $a \in \CC$, $|a| < 1$.  While rotations certainly preserve Euclidean
distance, the M\"{o}bius transformations do not---see Figure 1.

    \begin{figure}
    \centering
      \includegraphics[height=2.65in, width=2.75in]{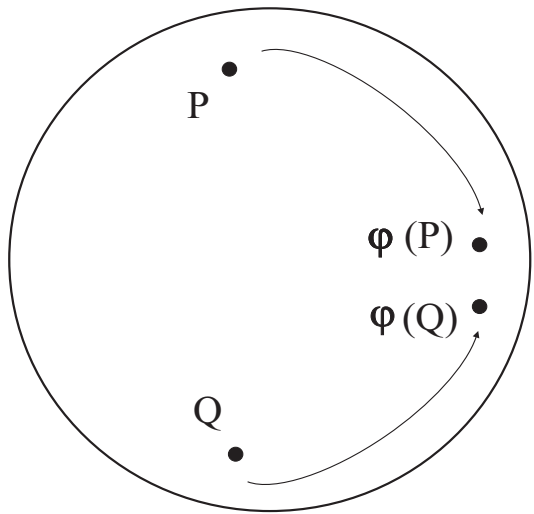}
      \caption{M\"{o}bius transformations do not preserve Euclidean distance.}
    \end{figure}

It is most convenient to describe the Poincar\'{e} metric in 
infinitesimal form.  In fact we set
$$
\rho(\zeta) = \frac{1}{1 - |\zeta|^2} \, .
$$
The {\it Poincar\'{e} length} of a vector $\xi$ at base point
$P$ is then defined to be
$$
\|\xi\|_{P, {\rm Poinc}} \equiv \rho(P) \cdot |\xi| \, ,
$$
where $|\xi|$ denotes Euclidean length of the vector $\xi$.  
Throughout most of this paper we shall use Finsler metric notation,
and we may as well begin now.  In that notation, we write the Poincar\'{e}
metric as
$$
F_{\cal P}^\Omega(P, \xi) = \frac{|\xi|}{1 - |P|^2} \, .
$$

Now we define the length of a piecewise $C^1$ curve $\gamma: [0,1] \rightarrow D$ to
be
$$
L_{\cal P}(\gamma) \equiv \int_0^1 F_{\cal P}^\Omega(\gamma(t), \gamma'(t)) \ dt \, .
$$
The {\it Poincar\'{e} distance} of two points $P$ and $Q$ in the disc, denoted $d_{\cal P}(P,Q)$, is
now declared to be the infinimum of the Poincar\'{e} lengths of all piecewise
$C^1$ curves connecting $P$ to $Q$.

The conformal invariance of the Poincar\'{e} metric is treated
in detail in the reference [KRA1]. We shall not discuss it
here. Suffice it to say that Poincar\'{e}'s construction is
special to the disc. It is in fact a matter of some interest
to equip virtually {\it any} domain in the plane (or in
higher-dimensional complex space) with a conformally or
biholomorphically invariant metric. And there are various
means of doing so. Certainly the most classical is to use the
uniformization theorem. We shall discuss that approach in the
next section. One of the first intrinsic techniques was
developed by Stefan Bergman in 1923 (see [BER]). Constantin
Carath\'{e}odory [CAR] created another in 1927. One of the
most recent is that developed in 1969 by S. Kobayashi [KOB1].
See also the definitive references [KOB2], [PFL], [KOB3],
[ROY].

Both Carath\'{e}odory's and Kobayashi's constructions have the advantage of
being elementary, intrinsic, flexible, and immediately accessible. Their
motivation from the proof of the Riemann mapping theorem is immediate. It
is a lovely example of modern mathematics at work. The present paper is
dedicated to the study of those two metrics.

It is a pleasure to thank the referees for very carefully reading
my manuscript and contributing considerable wisdom and insight.
The result is a cleaner and more precise presentation.

\section{The Uniformization Theorem}

The uniformization theorem of K\"{o}be and Poincar\'{e} is a remarkable generalization
of the Riemann mapping theorem.  We cannot prove it here, but we
provide a brief discussion.  

If $X$ is any topological space, then it has a simply connected universal
covering space $\widehat{X}$. The universal covering space is constructed
by fixing a point $x_0 \in X$ and considering the space of all paths in
$X$ emanating from $x_0$.  The covering map\index{universal covering space}
$$
\pi: \widehat{X} \ra X
$$
is a local homeomorphism.  We refer the reader to [SPA] or [HUS] for details.

In case $X$ is a domain $\Omega$ in the complex plane, 
or more generally a Riemann surface, then the universal
covering space $\widehat{X}$ will be a two-dimensional object (because
$\pi$ is a local homeomorphism), and $\widehat{X}$ can be endowed with a
complex structure by local pullback under $\pi$ of the complex structure
from $X$.  So $\widehat{X}$ is a simply connected analytic object.
What is it?

The uniformization theorem answers this question in a dramatic way.
Before we present the answer, let us first restate the question---stripped
of\index{universal covering space, complex structure on} 
all the preliminary material that led up to it.
\smallskip \\

\begin{quote}
{\bf QUESTION:}  What are all the simply connected Riemann
surfaces?
\end{quote}
The answer is
\begin{quote}
{\bf ANSWER:}  The only simply connected Riemann surfaces 
are {\bf (i)} the disk $D$, {\bf (ii)} the plane $\CC$,
and {\bf (iii)} the Riemann sphere $\widehat{\CC}$.
\end{quote}
\vspace*{.12in}

And\index{simply connected Riemann surfaces} 
in fact much more can be said.  Let us return to the
motivational discussion above.  If the original analytic
object $X$ is a sphere, then it turns out that 
the universal covering space $\widehat{X}$ will
be a sphere, and that is the {\it only} circumstance under
which\index{Riemann surfaces with the sphere as universal covering space} 
a sphere arises as the universal covering space.

If the original analytic object is a plane or a punctured
plane or a torus or a cylinder, then the universal covering
space $\widehat{X}$ is a plane, and these are the only
circumstances in which the plane arises as the universal
covering space.\index{Riemann surfaces with the plane as universal covering space}

In all other circumstances, the universal covering space is
the\index{Riemann surfaces with the disk as universal covering space} 
disk $D$.  In other words,
\begin{quote}
The universal covering space for any planar domain
except $\CC$ or $\CC \sm \{0\}$ is the disk $D$.
\end{quote}
This is powerful information, and those who study
Riemann surfaces have turned the result into
an important tool (see [FAK]).  

Suppose now that $U$ is a planar domain that is neither the entire plane
nor the punctured plane. Then the universal covering space is
(conformally) the disc and we have a covering map $\pi: D \ra U$. Then we
may push the Poincar\'{e} metric from the disc down to $U$---that is to
say, measure the length of a tangent vector to $U$ at $P \in U$ by pulling
the vector back up to $D$ by way of $\pi$. And so virtually {\it any}
planar domain may be equipped with an invariant metric. We call such a
domain {\it hyperbolic}---see [KRA1]. One of the points of the present
paper is that there are instead distinct advantages to constructing the
invariant metric intrinsically.  The constructions of Carath\'{e}odory
and Kobayashi in fact generalize to a broad range of circumstances---even
complex manifolds of many variables---and have proved to be powerful tools
for function theory.

\section{Motivation by Way of the Schwarz and Schwarz-Pick Lemmas}

The construction of the Carath\'{e}odory and Kobayashi metrics
is motivated in a natural way by Schwarz's lemma.  The fact that
an invariant metric is thereby constructed is closely related to
the more general Schwarz-Pick lemma.  We take this opportunity to
review those ideas.

The classical Schwarz lemma is part of the grist of every
complex analysis class.  A version of it says this:

\begin{lemma} \sl
Let $f: D \ra D$ be holomorphic.  Assume that $f(0) = 0$.
Then
\begin{enumerate}
\item[{\bf (a)}]  $|f(z)| \leq |z|$ for all $z \in D$;
\item[{\bf (b)}]  $|f'(0)| \leq 1$.
\end{enumerate}
At least as important as these two statements are
the cognate {\it uniqueness statements}:
\begin{enumerate}
\item[{\bf (c)}]  If $|f(z)| = |z|$ for some $z \ne 0$, then
$f$ is a rotation:  $f(z) = \lambda z$ \hfill \break
\null \indent for some unimodular complex constant $\lambda$;
\item[{\bf (d)}]  If $|f'(0)| = 1$, then $f$ is a rotation:  $f(z) = 
\lambda z$ for some \hfill \break
\null \indent unimodular complex constant $\lambda$.
\end{enumerate}
\end{lemma}

There\index{Schwarz lemma, classical form} 
are a number of ways to prove this result.  The classical
argument is to consider $g(z) = f(z)/z$.  On a circle
$|z| = 1 - \e$, we see that $|g(z)| \leq 1/(1 - \e)$.  Thus
$|f(z)| \leq |z|/(1 - \e)$.  Since this inequality holds
for all $\e > 0$, part {\bf (a)} follows.  
The Cauchy estimates show that $|f'(0)| \leq 1$.

For the uniqueness, if $|f(z)| = |z|$ for some $z \ne 0$,
then $|g(z)| = 1$.  The maximum modulus principle then
forces $|f(z)| = |z|$ for all $z$, and hence $f$ is a rotation.
If instead $|f'(0)| = 1$, then $|g(0)| = 1$ and again
the maximum modulus principle yields that $f$ is a rotation.

The Schwarz-Pick lemma observes that there is no need
to restrict to $f(0) = 0$.  Once one comes up with the 
right\index{Schwarz-Pick lemma} formulation, the proof is straightforward:

\begin{proposition} \sl
Let $f:D \ra D$.  Assume that $a \ne b$ are elements of $D$
and that $f(a) = \alpha$, $f(b) = \beta$.  Then
\begin{enumerate}
\item[{\bf (a)}]  $\displaystyle 
 \left | \frac{\beta - \alpha}{1 - \overline{\alpha}\beta} \right |
    \leq \left | \frac{b - a}{1 - \overline{a} b} \right |$;
\item[{\bf (b)}]  $\displaystyle
 \left | f'(a) \right | \leq \frac{1 - |\alpha|^2}{1 - |a|^2}$.
\end{enumerate}
There is also a pair of uniqueness statements:
\begin{enumerate}
\item[{\bf (c)}]  If $\displaystyle 
 \left | \frac{\beta - \alpha}{1 - \overline{\alpha}\beta} \right |
    = \left | \frac{b - a}{1 - \overline{a} b} \right |$,
then $f$ is a conformal self-map of the disk $D$;
\item[{\bf (d)}]  If $\displaystyle
 \left | f'(a) \right | = \frac{1 - |\alpha|^2}{1 - |a|^2}$,
then $f$ is a conformal self-map of the disk $D$.
\end{enumerate}
\end{proposition}
{\bf Proof:} We sketch the proof.  Recall that, for
$a$ a complex number in $D$,
$$
\varphi_a(\zeta) = \frac{\zeta - a}{1 - \overline{a}\zeta}
$$
defines a {\small \it M\"{o}bius transformation}.  This is a conformal
self-map of the disk that takes $a$ to 0.  Note that $\varphi_{-a}$
is the inverse mapping to $\varphi_a$.

Now, for the given $f$, consider
$$
g(z) = \varphi_\alpha \circ f \circ \varphi_{-a} \, .
$$
Then $g: D \ra D$ and $g(0) = 0$.  So the standard Schwarz
lemma applies to $g$.  By part {\bf (a)} of that lemma,
$$
|g(z)| \leq |z| \, .
$$
Letting $z = \varphi_a(\zeta)$ yields
$$
|\varphi_\alpha \circ f(\zeta)| \leq |\varphi_a(\zeta)| \, .
$$
Writing this out, and setting $\zeta = b$, gives the conclusion
$$
\left | \frac{\beta - \alpha}{1 - \overline{\alpha}\beta} \right |
    \leq \left | \frac{b - a}{1 - \overline{a} b} \right | \, .
$$
That is part {\bf (a)}.

For part {\bf (b)}, we certainly have that
$$
\left | (\varphi_\alpha \circ f \circ \varphi_{-a} )'(0) \right | \leq 1 \, .
$$
Using the chain rule, we may rewrite this as
$$
\left | \varphi'_\alpha(f \circ \varphi_{-a}(0)) \right | \cdot
   \left | f'(\varphi_{-a}(0)) \right | \cdot
  \left | \varphi'_{-a} (0) \right | \leq 1 \, .  \eqno (*)
$$
Now of course 
$$
\varphi'_a(\zeta) = \frac{1 - |a|^2}{(1 - \overline{a}\zeta)^2} \, .
$$
So we may rewrite $(*)$ as
$$
\left ( \frac{1 - |\alpha|^2}{(1 - |\alpha|^2)^2} \right ) \cdot 
     |f'(a)| \cdot (1 - |a|^2) \leq 1 \, .
$$
Now part {\bf (b)} follows.

We leave parts {\bf (c)} and {\bf (d)} as exercises for the reader.
\endpf
\smallskip  \\

\rm
The quantity
$$
\rho(a,b) = \frac{|a - b|}{|1 - \overline{a}b|}
$$
is called the {\it pseudohyperbolic metric}.  It is actually
a metric on $D$ (details left to the reader).  It is not
identical to the Poincar\'{e}-Bergman metric.  In fact
it is not a Riemannian metric\index{pseudohyperbolic metric} 
at all.  But it is still
true that conformal maps of the disk are distance-preserving
in the pseudohyperbolic metric.  Exercise:  Use the Schwarz-Pick
lemma to prove this last assertion.

One useful interpretation of the Schwarz-Pick lemma is that
a holomorphic function $f$ from the disk to the disk must take
each disk $D(0,r)$, $0 < r < 1$, into (but not necessarily onto)
the image of that disk under the linear fractional map
$$
z \mapsto \frac{z + \alpha}{1 + \overline{\alpha}z} \, ,
$$
where $f(0) = \alpha$.  This image is in fact (in case $-1 < \alpha < 1$) a standard
Euclidean disk with center on the real axis at $\alpha$ and
diameter (in case $0 < \alpha < 1$) given by the interval
$$
\left [ \frac{\alpha - r}{1 - \alpha r} \ , \
           \frac{\alpha + r}{1 + \alpha r} \right ] \, .
$$

The reader will see, when encountering the definitions of 
the Carath\'{e}odory and Kobayashi metrics, the Schwarz lemma
acting as motivation.  Certainly the Schwarz lemma arises frequently
in the proofs of the basic results about these metrics.

\section{Basic Facts about the Kobayashi Metric}

Following the paradigm set in Section 0 (for the Poincar\'{e} metric), we shall define		
the Kobayashi metric at first on the infinitesimal level.  That
is to say, we shall specify the length of a tangent vector at each point.
We will always let $\Omega$ denote a connected, open set, or a {\it domain}.
Following tradition, we let $D = \{\zeta \in \CC: |\zeta| < 1\}$ denote
the unit disc and we let $\Omega(D)$ denote the collection of
holomorphic functions from $D$ to $\Omega$.  If $z \in \Omega$ then we further
let $\Omega^z(D)$ denote the subcollection of elements $f$ of $\Omega(D)$ which
satisfy $f(0) = z$. 

So let $\Omega \subseteq \CC$ be a domain.
Fix a point $P \in \Omega$ and a vector $\xi$ which is thought of as
being tangent to the plane at the point $P$.  We define
the infinitesimal Kobayashi or Kobayashi/Royden length of $\xi$ at $P$ to be
\begin{eqnarray*}
F_K^\Omega(P,\xi) & \equiv & \inf\{\alpha: \alpha > 0\ \mbox{\rm and} \ \exists 
f \in \Omega(D)\ \mbox{\rm with} \ f(0) = P,  f'(0) = \xi/\alpha \} \\
    & = & \inf \left\{\frac{|\xi|}{|f'(0)|}: f \in \Omega^P (D) \right \} \, .
\end{eqnarray*}
It is in general not the case that $F_K^\Omega$ satisfies a triangle inequality
in the second entry.  Nonetheless we can, as indicated in our discussion
of the Poincar\'{e} metric, construct from it a useful metric.
\smallskip \\

\noindent {\bf Remark:}  Recall the standard, modern proof of the Riemann mapping theorem (see [GRK]).
We are given a simply connected domain $\Omega$ (not all of $\CC$), and our goal
is to construct a conformal mapping of $D$ to $\Omega$.  We fix a point $P \in \Omega$ and we 
consider the family ${\cal S}$ of holomorphic mappings $\varphi: D \ra \Omega$ with $\varphi(0) = P$.
A normal families argument is used to show that there is an element $\varphi^*$ of ${\cal S}$ that
maximizes the modulus of the derivative at 0.  The function $\varphi^*$ turns out to
be the conformal mapping that we seek.

Now look at the definition of the Kobayashi/Royden metric.  The metric at a point $P$ in
the direction $\xi$ minimizes the expression $|\xi|/|f'(0)|$ over mappings
$f$ of the disc into $\Omega$.  This is the same as {\it maximizing} the quantity
$|f'(0)|$.  Thus we see the proof of the Riemann mapping theorem coming back to life
in the definition of the Kobayashi/Royden metric.\footnote{In fact the idea behind the
Kobayashi metric has a long history.  Even in the 1920s, T. Rado observed that the
same extremal problem may be used to produce a proof of the uniformization theorem
for planar domains.  See [GOL, p.\ 256].}
\endpf
\smallskip \\

\begin{definition}    \rm
Let $\Omega \subseteq \CC$ be open and $\gamma: [0,1] \ra \Omega$ a piecewise $C^1$ curve.
The {\it Kobayashi/Royden length} of $\gamma$ is defined to be\footnote{A word
needs to be said about why $F_K^\Omega$ is integrable.  In fact it is not difficult
to see that $F_K^\Omega$ is lower semicontinuous, since it is the infimum of
continuous functions.  And that is sufficient for the integrability.}
$$ 
L_K^\Omega(\gamma) = L_K(\gamma) = \int_0^1 F_K^\Omega(\gamma(t), \gamma'(t)) dt . 
$$
\end{definition}

\begin{definition}  \rm
Let $\Omega \subseteq \CC$ be an open set and $z,w \in \Omega.$  The
{\it (integrated) Kobayashi/Royden distance} between $z$ and $w$ is defined to be
$$  
K^\Omega(z,w) = K(z,w) = \inf\{L_K(\gamma): \gamma\ \mbox{is a piecewise $C^1$ curve connecting} \
                    \ z\ \mbox{and} \ w\} . 
$$
\end{definition}

Of course it must be noted that $K^\Omega$ need not be a distance function in the
classical sense of the term.  As an instance, if $\Omega$ is the entire complex
plane then $K^\Omega$ is identically equal to zero.  A domain for which
$K^\Omega$ is a genuine (nondegenerate) distance is called {\it hyperbolic} (this
is equivalent to our earlier use of the term ``hyperbolic'').
The book [KRA1] has a concise treatment of hyperbolicity.  See also [KOB2].  For all practical
purposes, hyperbolicity is an invariant version of boundedness---in other words,
a hyperbolic domain has all the key properties of a bounded domain, but hyperbolicity
has the additional advantage of being invariant under conformal mappings.

One of the most important, and most interesting, properties of
the Kobayashi metric is that a holomorphic function is
distance decreasing in the metric. We shall be able to make
good use of it in the examples below.

\begin{proposition}[The Distance Decreasing Property of the Kobayashi Metric]    
If $\Omega_1, \Omega_2$ are domains in $\CC, z, w \in \Omega_1, \xi \in \CC,$ and if $f: \Omega_1 \ra \Omega_2$ is holomorphic,
then
$$
F_K^{\Omega_2}(f(z),f'(z)\xi) \leq F_K^{\Omega_1}(z,\xi) \quad \hbox{and} \quad
   K^{\Omega_2}(f(z), f(w)) \leq K^{\Omega_1}(z,w) \, .
$$
\end{proposition}

\noindent {\bf Remark:}  Observe that the Chain Rule demands
that we put a factor of $f'(z)$ in front of the
tangent vector when we calculate $F_K^{\Omega_2}(f(z), \, \cdot \, )$.
\endpf
\smallskip \\

\noindent {\bf Proof of the Proposition:}  We prove the first inequality and leave the second for the reader.

Let $\varphi: D \ra \Omega_1$ satisfy $\varphi(0) = z$.  We call $\varphi$ a {\it candidate
mapping} for the Kobayashi metric at the point $z$ on the domain $\Omega_1$.  Then $f \circ \varphi$
is a candidate mapping for the Kobayashi metric at the point $f(z)$ on the domain $\Omega_2$.
Thus
$$
F_K^{\Omega_2}(f(z), f'(z)\xi) = \inf_{g \in \Omega^{f(z)}_2(D)} \frac{|f'(z)\xi|}{|g'(0)|} \leq
      \frac{|f'(z)\xi|}{|(f \circ \varphi)'(0)|} = \frac{|\xi|}{|\varphi'(0)|} \, .
$$
Now we take the infimum over all candidates $\varphi$ to obtain
$$
F_K^{\Omega_2}(f(z), f'(z)\xi) \leq F_K^{\Omega_1}(z, \xi) \, .  \eqno \BoxOpTwo
$$

\begin{corollary}     
If $f: \Omega_1 \ra \Omega_2$ is conformal then $f$ is an isometry in 
the Kobayashi/Royden metric.  This means that
$f$ preserves distances:
$$
F_K^{\Omega_1}(z, \xi) = F_K^{\Omega_2}(f(z), f'(z)\xi) \quad \hbox{and} \quad K^{\Omega_2}(f(z), f(w)) = K^{\Omega_1} (z,w) \, .
$$
\end{corollary}

\noindent {\bf Remark:}  A caution is in order here.  The reader who knows some differential
geometry will be accustomed to the term ``isometry'', and will think of such a mapping
as preserving distances in a strong (classical) sense.  The ``metrics'' that we consider
now may degenerate to 0, so our present use of the term ``isometry'' is somewhat more general.
\endpf
\smallskip \\

\noindent {\bf Proof of the Corollary:}  Let us prove the first assertion.  We leave the second to the reader.  Now the
proposition certainly tells us that
$$
F_K^{\Omega_2}(f(z), f'(z)\xi) \leq F_K^{\Omega_1}(z,\xi) \, .  \eqno (*)
$$
But we may also apply the proposition to $f^{-1}: \Omega_2 \ra \Omega_1$.  The result
is 
$$
F_K^{\Omega_1}(f^{-1}(a), [f^{-1}]'(a) \tau) \leq F_K^{\Omega_2}(a,\tau) \, .
$$
Now simply let $a = f(z)$ and $\tau = f'(f^{-1}(a))\xi$ to obtain
$$
F_K^{\Omega_1}(z, \xi) \leq F_K^{\Omega_2}(f(z),f'(z)\xi) \, .  \eqno (**)
$$
Combining $(*)$ and $(**)$ yields
$$
F_K^{\Omega_1}(z, \xi) = F_K^{\Omega_2}(f(z),f'(z)\xi) \, .  \eqno \BoxOpTwo
$$

\begin{corollary}    
If $\Omega_1 \subseteq \Omega_2 \subseteq \CC$ then, for any $z,w \in \Omega_1,$
any $\xi \in \CC,$ we have
$$
F_K^{\Omega_1}(z,\xi) \geq F_K^{\Omega_2}(z,\xi) \qquad \hbox{and} \qquad
     K^{\Omega_1}(z,w) \geq K^{\Omega_2}(z, w) .
$$
\end{corollary}
{\bf Proof:}  Simply apply the proposition to the inclusion
mapping $i: \Omega_1 \ra \Omega_2$.
\endpf
\smallskip \\

\section{Basic Facts about the \cara\ Metric}

Following the model set in the last section, we shall define the
\cara\ metric at first on the infinitesimal level. That is to say, we
shall specify the length of a tangent vector at each point. As usual, we
let $\Omega$ denote a connected, open set, or a domain. Following
tradition, we let $D(\Omega)$ denote the collection of holomorphic
functions\footnote{Of course it is possible that $D(\Omega)$ is trivial---for
example when $\Omega$ is the entire plane.} from $\Omega$ to $D$. If $z \in \Omega$ then we further let
$D^z(\Omega)$ denote the subcollection of elements $g$ of $D(\Omega)$ such that
$g(z) = 0$.

So let $\Omega \subseteq \CC$ be a domain.
Fix a point $P \in \Omega$ and a vector $\xi$ which is thought of as
being tangent to the plane at the point $P$.  We define
the infinitesimal \cara\ length of $\xi$ at $P$ to be
$$
F_C^\Omega(P,\xi) \equiv \sup_{f \in D^P (\Omega) 
                       \atop
                       f(P) = 0}
                                  |f'(P) \xi| \, .
$$

\noindent {\bf Remark:}  Refer to the discussion in the Remark following
the definition at the beginning of Section 3 of the Kobayashi metric.
It is worth noting that we could as well prove the Riemann mapping theorem
by considering maps of the domain $\Omega$ into the disc $D$ and maximizing the
derivative at the point $P$.
Now look at the definition of the \cara\ metric.  The metric at a point $P$ in
the direction $\xi$ maximizes the expression $|f'(P) \xi|$ over mappings
$f$ of $\Omega$ into the disc.  This is the same as {\it maximizing} the quantity
$|f'(P)|$.  Thus we see the proof of the Riemann mapping theorem coming back to life
in the definition of the \cara\ metric.

It is worth noting that an extremal function for the Carath\'{e}odory metric
always exists, as can be seen with a normal families argument.  The extremal
function is often termed the {\it Ahlfors function}.   It is, in many respects,
a generalization of the Riemann mapping function (which is what it is in
case $\Omega$ is simply connected).  See [FIS] or [KRA3] for a consideration of
the Ahlfors function.
\endpf
\smallskip \\

\begin{definition}    \rm
Let $\Omega \subseteq \CC$ be open and $\gamma: [0,1] \ra \Omega$ a piecewise $C^1$ curve.
The {\it \cara\ length} of $\gamma$ is defined to be
$$ 
L_C^\Omega(\gamma) = L_C(\gamma) = \int_0^1 F_C^\Omega(\gamma(t), \gamma'(t)) dt . 
$$
\end{definition}
We note that $F_C$ is integrable for reasons similar to the ones
given for $F_K$ in Section 1.

Next we are going to define the integrated \cara\ distance in $\Omega$.  But
now our approach will not parallel that for the Kobayashi metric.  In fact
we want the \cara\ metric to have a certain ``minimal property'' among
all metrics for which holomorphic functions are distance decreasing.  This
necessitates a new approach.

\begin{definition}  \rm
Let $\Omega \subseteq \CC$ be an open set and $z,w \in \Omega.$   The
{\em Carath\'{e}odory distance} between $z$ and $w$ is defined to be
$$
  C^\Omega(z,w) = \sup_{f \in D(\Omega)} d_{\cal P}(f(z), f(w)) , 
$$
where $d_{\cal P}$ is the Poincar\'{e} distance on $D$.
\end{definition}

\noindent {\bf Remark:}  Of course the Carath\'{e}odory distance can
be trivial---for instance if $\Omega$ is the entire plane.
\endpf
\smallskip \\

One of the most important, and most interesting, properties of the \cara\ metric
is that a holomorphic function is distance decreasing in the metric.  We shall
be able to make good use of it in the examples below.

\begin{proposition}[The Distance Decreasing Property of the \cara\ Metric]    
If $\Omega_1, \Omega_2$ are domains in $\CC, z, w \in \Omega_1, \xi \in \CC,$ and if $f: \Omega_1 \ra \Omega_2$ is holomorphic,
then
$$
F_C^{\Omega_2}(f(z),f'(z)\xi) \leq F_C^{\Omega_1}(z,\xi) \quad \hbox{and} \quad
   C^{\Omega_2}(f(z), f(w)) \leq C^{\Omega_1}(z,w) \, .
$$
\end{proposition}

\noindent {\bf Remark:}  Observe that the Chain Rule demands
that we put a factor of $f'(z)$ in front of the
tangent vector when we calculate $F_C^{\Omega_2}(f(z), \, \cdot \, )$.
\endpf
\smallskip \\

\noindent {\bf Proof of the Proposition:}  We prove the first inequality and leave the second for the reader.

Let $\varphi: \Omega_2 \ra D$ satisfy $\varphi(f(z)) = 0$.  We call $\varphi$ a {\it candidate
mapping} for the \cara\ metric at the point $f(z)$ on the domain $\Omega_2$.  Then $\varphi \circ f$
is a candidate mapping for the \cara\ metric at the point $z$ on the domain $\Omega_1$.
Thus
$$
F_C^{\Omega_1}(z, \xi) = \sup_{g \in D^z(\Omega_1)} |g'(z)\xi| \geq
      |(\varphi \circ f)'(z) \xi| = |\varphi'(0)| \cdot |f'(z)| \cdot |\xi| \, .
$$
Now we take the supremum over all candidates $\varphi$ to obtain
$$
F_C^{\Omega_1}(z, \xi) \geq F_C^{\Omega_2}(f(z), f'(z) \xi) \, . \eqno \BoxOpTwo
$$

\begin{corollary}     
If $f: \Omega_1 \ra \Omega_2$ is conformal then $f$ is an isometry in the
Carath\'{e}odory metric.  This means that
$f$ preserves distances:
$$
F_C^{\Omega_1}(z, \xi) = F_C^{\Omega_2}(f(z), f'(z)\xi) \quad \hbox{and} \quad C^{\Omega_2}(f(z), f(w)) = C^{\Omega_1} (z,w) \, .
$$
\end{corollary}
{\bf Proof:}  Let us prove the second assertion.  We leave the first to the reader.  Now the
proposition certainly tells us that
$$
C^{\Omega_2}(f(z), f(w)) \leq C^{\Omega_1}(z,w) \, .  \eqno (*)
$$
But we may also apply the proposition to $f^{-1}: \Omega_2 \ra \Omega_1$.  The result
is 
$$
C^{\Omega_1}(f^{-1}(a), f^{-1}(b)) \leq C^{\Omega_2}(a,b) \, .
$$
Now simply let $a = f(z)$ and $b = f(w)$ to obtain
$$
C^{\Omega_1}(z, w) \leq C^{\Omega_2}(f(z),f(w)) \, .  \eqno (**)
$$
Combining $(*)$ and $(**)$ yields
$$
C^{\Omega_1}(z, w) = C^{\Omega_2}(f(z),f(w)) \, .   \eqno \BoxOpTwo
$$

\begin{corollary}    
If $\Omega_1 \subseteq \Omega_2 \subseteq \CC$ then for any $z,w \in \Omega_1,$
any $\xi \in \CC,$ we have
$$
F_C^{\Omega_1}(z,\xi) \geq F_C^{\Omega_2}(z,\xi) \qquad \hbox{and} \qquad
     C^{\Omega_1}(z,w) \geq C^{\Omega_2}(z, w) .
$$
\end{corollary}
{\bf Proof:}  Simply apply the proposition to the inclusion
mapping $i: \Omega_1 \ra \Omega_2$.
\endpf
\smallskip \\

\section{Comparison of the Kobayashi and \cara\ Metrics}

First, it is always the case that the Kobayashi metric majorizes the
\cara\ metric:

\begin{proposition} \sl
Let $\Omega \subseteq \CC$ be a domain.  Let $P \in \Omega$ and let $\xi$ be a vector.
Then
$$
F_C^\Omega(P, \xi) \leq F_K^\Omega(P,\xi) \, .
$$
\end{proposition}
{\bf Proof:}  Let $\varphi: D \ra \Omega$ be a candidate mapping for the Kobayashi
 metric at $P \in \Omega$.  Let $\psi: \Omega \ra D$ be a candidate mapping
 for the \cara\ metric at $P \in \Omega$.  Then $h \equiv \psi \circ \varphi:  D \ra D$
 and $h(0) = 0$.  By Schwarz's lemma, $|h'(0)| \leq 1$.  But this just says that
 $$
 |\psi'(P)| \leq \frac{1}{|\varphi'(0)|} \, .
 $$
 Now first take the infimum on the right over all candidate functions $\varphi$ for
 the Kobayashi metric then take the supremum on the left over all candidate functions
 $\psi$ for the \cara\ metric.  The result is
$$
F_C^\Omega(P, \xi) \leq F_K^\Omega(P,\xi) \, .   \eqno \BoxOpTwo
$$

We conclude the section with an interesting extremal property of $C^\Omega$.

\begin{theorem}  \sl
Let $\Omega \subseteq \CC$ be an open set.  Let $d$ be any metric on $\Omega$
that satisfies $d(z,w) \geq d_{\cal P}(f(z),f(w))$ for all $f \in D(\Omega)$ and
$z,w \in \Omega.$  Then $d(z,w) \geq C^\Omega(z,w).$
\end{theorem}
{\bf Proof:}  Exercise.  Use the definition of $C^\Omega$.
\endpf
\smallskip \\

It is worth noting that the Kobayashi metric satisfies an analogous
extremal property:

\begin{theorem} \sl
Let $\Omega \subseteq \CC$ be an open set.  Let $d$ be any metric on $\Omega$
that satisfies $d(z,w) \geq d_{\cal P}(f(z),f(w))$ for all $f \in D(\Omega)$ and
$z,w \in \Omega.$  Then $d(z,w) \leq K^\Omega(z,w).$
\end{theorem}

We leave the details of this last result for the interested reader,
or consult [KRA2].

In the present paper the roles of the Carath\'{e}odory and Kobayashi
metrics are virtually interchangeable. Any proof that uses the
Carath\'{e}odory metric could just as well use the Kobayashi metric, and
vice versa. But both metrics are interesting because they are defined in a
dual manner, and because the one (the Kobayashi) always majorizes the
other (the Carath\'{e}odory).  As we have already noted, the
Kobayashi metric is the {\it largest} metric in which holomorphic mappings
are distance decreasing and the Carath\'{e}odory metric is the smallest.
It is an interesting, and more recent, result of K. T. Hahn (see [KRA2])
that the Berman metric always majorizes the Carath\'{e}odory metric.  It is
also known, thanks to an example of Diederich and Forn\ae ss [DIF], that the Bergman
metric {\it cannot} in general be compared with the Kobayashi metric.

\section{Calculation of the \cara\ and Kobayashi Metrics}

Precious little is known about explicitly calculating the Kobayashi/Royden
metric or the \cara\ metric. For special domains such as the disc or the
annulus, the automorphism group is a powerful tool for obtaining explicit
formulas.\footnote{Of course we
{\it could} write down a formula for the Kobayashi or \cara\
metric of the upper halfplane, for example. But that is only
because the halfplane is---by way of the Cayley
map---conformally equivalent to the unit disc. However there
is an extensive literature on this subject. The paper [SIM]
considers the Carath\'{e}odory metric on the annulus in some
detail. Gehring and Palka [GEP] have developed a {\it
quasihyperbolic metric} that can be used, with comparison
arguments, to obtain estimates for the Kobayashi metric.}  In many circumstances one can instead {\it estimate} the
metrics, and that is sufficient for applications (see, for example,
[GEP]). Let us now, just for illustrative purposes, calculate the Kobayashi metric on the
disc.

\begin{example} \rm
We let $\Omega = D$ be the unit disc.  We begin by calculating the infinitesimal
Kobayashi metric at the origin 0.  Let $\xi = 1 + i0$.  We calculate 
$F_K^D(0,\xi)$.  So let 
$$
f: D \ra \Omega
$$
be holomorphic and satisfy $f(0) = 0$.  By Schwarz's lemma, we know
that $|f'(0)| \leq 1$.  But in fact the function $f_0(\zeta) \equiv \zeta$
maps $D$ to $\Omega$ with $f(0) = 0$ and $f'_0(0) = 1$.  We conclude that 1 is the
extremal value and
$$
F_K^D(0, \xi) = 1 \, .
$$

Now if $\xi$ is any vector of length 1 then $\xi = e^{i\theta}$ for some
$0 \leq \theta < 2\pi$.  Then
$$
F_K^D (0, \xi) = F_K^D(0, e^{i\theta} \cdot 1) = |e^{i\theta}| F_K^D (0, 1) = 1 \, .
$$
Next, if $\xi$ is {\it any} vector then write $\xi = r e^{i\theta}$.
Thus
$$
F_K^D(0,\xi) = F_K^D(0, r e^{i\theta}) = r F_K^D (0, e^{i\theta}) = r = |\xi| \, .
$$
[Again we remind the reader that $|\xi|$ denotes the {\it Euclidean} length
of the vector $\xi$.]

Our next task is to derive a formula for $F_K^D$ at an arbitrary base point
$P \in D$.  Now notice that the M\"{o}bius transformation
$$
\varphi(\zeta) = \frac{\zeta - P}{1 - \overline{P}\zeta}
$$ 
maps $P$ to $0$.  Also
$$
\varphi'(P) = \frac{(1 - \overline{P} \zeta)\cdot 1 - (\zeta - P) \cdot (-\overline{P})}{(1 - \overline{P}\zeta)^2} \biggr |_{\zeta = P}
   = \frac{1}{1 - |P|^2} \, .
$$
Therefore we may calculate, for any vector $\xi$, that
$$
F_K^D(P, \xi) = F_K^D(\varphi(P), \varphi'(P) \xi) = F_K^D(0, \xi/(1 - |P|^2))
    = \frac{1}{1 - |P|^2} \cdot F_K^D(0, \xi) = \frac{|\xi|}{1 - |P|^2} \, .   \eqno \BoxOpTwo
$$
\end{example}

Of course what we have just calculated is the {\it infinitesimal form} of the
Kobayashi metric.  It is certainly of interest to have a formula for the integrated
form---as that will be a genuine metric in the classical sense.  And it will
be invariant under conformal mappings.  We note already that the conclusion
of Example 5 already shows that, on the disc, the Kobayashi metric coincides with
the Poincar\'{e} metric.  And in fact a calculation nearly identical to the one
we just performed shows that the \cara\ metric coincides with the Kobayashi and
Poincar\'{e} metrics on the unit disc.

\begin{proposition} \sl
The length of the curve $\gamma(t) = \gamma_\epsilon(t) = t + i0$, $0 \leq t \leq 1 - \epsilon$ in the
Kobayashi metric on the the disc $D$ is
$$
L_K^D(\gamma_\epsilon) = \frac{1}{2} \cdot \log \left [ \frac{2 - \epsilon}{\epsilon} \right ] \, .
$$
\end{proposition}

\noindent {\bf Remark:}  This proposition is particularly interesting, for it tells
us that
$$
\lim_{\epsilon \ra 0^+} L_K^D(\gamma_\epsilon) = + \infty \, .
$$
In other words, the distance from the origin to the boundary of $D$---at least along a straight
line segment---is $+\infty$.
If we can show that the straight line segment is the shortest curve in
the Kobayashi metric from $0$ to $(1 - \epsilon) + i0$, then we will have proved that
the distance of 0 to $\partial D$ is $+\infty$.  [We shall establish this latter contention
in a moment.  First we prove the proposition.]	 This will, in turn, say that
the unit disc $D$ is {\it complete} in the Kobayashi metric.

At first such a statement may seem bewildering.  How can a bounded, open set be complete?  It certainly
does not appear to be closed in any sense; on the contrary, it is open!  But think of the Euclidean
plane in the ordinary Euclidean metric.  It is certainly complete.  And that is because the 
boundary is infinitely far away.  That is exactly what is happening with the Kobayashi metric
on the unit disc.
\endpf
\smallskip \\

\noindent {\bf Proof of the Proposition:}  

Now
\begin{eqnarray*}
L_K^D(\gamma)
& = & \int _0^{1 - \epsilon} F_K^D(\gamma(t), \gamma'(t)) \, dt  \\
& = & \int_0^{1 - \epsilon} {| \gamma' (t) | \over 1 - | \gamma
(t) | ^2} dt \\
& = & \int _0^{1 - \epsilon} {1 \over 1 - t^2}  dt \\
& = & \frac{1}{2}  \log \left [ {2 - \epsilon \over \epsilon} \right ] \, . \qquad \qquad \BoxOpTwo
\end{eqnarray*}

\begin{proposition} \sl  Among all continuously differentiable curves of
the form
$$ 
\mu (t) = t + iw(t) \ ,\qquad \ 0  \leq   t \leq   1 - \epsilon \ , 
$$ 
\noindent that satisfy  $\mu (0) = 0$  and  $\mu (1- \epsilon)
= 1 - \epsilon + 0i,$  the one of least length in the Kobayashi metric is  $\gamma (t) = t.$
Here  $w(t)$ is any continuously differentiable, real-valued
function.
\end{proposition}
{\bf Proof:}  In fact, for any such candidate $\mu$, we have
\begin{eqnarray*}
L_K^D (\mu )
& = & \int_0^{1 - \epsilon} 
F_K^D(\mu(t), \mu'(t)) \, dt \\
& = &  \int_0^{1 - \epsilon} {1 \over 1 - | \mu (t) |^2} 
\cdot | {\dot \mu} (t) |  \, dt  \\
& = &  \int_0^{1 - \epsilon} {1 \over 1 - t^2 - [w(t)]^2}
\cdot (1 + [w' (t)]^2)^{1/2} \, dt \, . 
\end{eqnarray*}
\noindent However
$$
{1 \over 1 - t^2 - [w(t)]^2} \geq  {1 \over 1 - t^2} 
\ \ \hbox{and}\ \  (1 + [w'(t)]^2)^{1/2} \geq 1.
$$
\noindent We conclude that
$$
L_K^D (\mu) \geq
\int_0^{1 - \epsilon} {1 \over 1 - t^2} dt
= L_K^D (\gamma).
$$
\noindent This is the desired result.

Notice that, with only small modifications, this argument can also
be applied to  {\it piecewise\/}   continuously differentiable
curves $t + iw(t).$  
\endpf
\smallskip  \\

In fact if a piecewise continuously differentiable curve connecting
the point  $0 \in D$  to  $(1 - \epsilon) + 0i \in  D$  is {\it
not\/} of the
form
$$
\mu (t) = t + iw(t), \eqno (*)
$$
\noindent then it may cross itself.  Of course
we can eliminate the loops and thereby create a shorter curve.
If the resulting curve is still not the graph of a function,
then elementary comparisons show that it 
will be longer than a curve of the form
$(*)$  (see Figure 2).  We may conclude that the curve $\gamma$ 
in the proposition is the shortest of all curves connecting  $0$  to 
$(1 - \epsilon) + 0i.$

    \begin{figure}
    \centering
      \includegraphics[height=2.65in, width=2.75in]{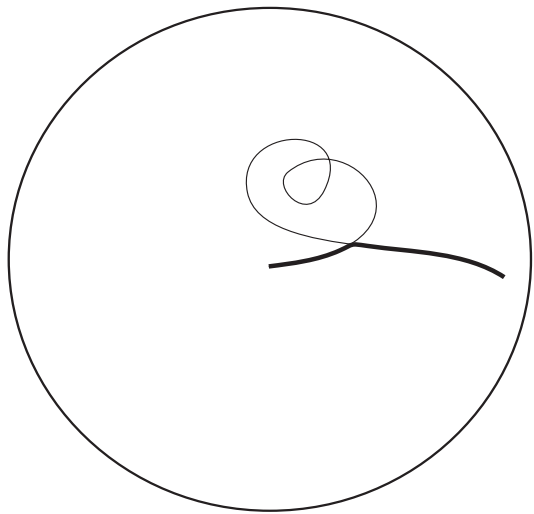}
      \caption{A curve that is not a graph.}
    \end{figure}

Of course we can use the last result to give an explicit formula
for the Kobayashi or \cara\ metric on the disc.  This we now do.

\begin{proposition} \sl
The integrated Kobayashi or \cara\ distance of two points
$P$ and $Q$ in $D$ is given by
$$
d(P,Q) = {1 \over 2} \log
\left ( {{1 + \bigl |{{P - Q} \over 1 - \overline{P}Q}\bigr |} 
\over {1 - \bigl |{{P - Q} \over {1 - {\overline{P}Q}} }\bigr | }} \right )
 \, . 
$$
\end{proposition}
{\bf Proof:}  In case  $P = 0$  and  $Q = R + i0$,  the result was 
already noted in Proposition 9.
In the general case, note that we may define
$$
\varphi (z) = {z - P \over 1 - \overline{P}z} \ ,
$$
\noindent a M\"obius transformation of the disc.
Then, by Corollaries 4 and 7 (letting $d$ denote either
the Kobayashi or the \cara\ distance),
$$
d (P,Q) = d  (\varphi (P), \varphi (Q)) 
= d (0,\varphi (Q)).
$$
\noindent Next we have
$$
d (0,\varphi (Q)) = d (0, | \varphi (Q) | ) \eqno
(*)
$$
\noindent since there is a rotation of the disc taking  
$\varphi (Q)$ to $|\varphi (Q)| + i0$.
Finally,
$$
| \varphi (Q) | 
= \left | {P - Q \over 1 - \overline{P}Q} \right | \ ,
$$
\noindent so that  $(*)$  together with the special case treated
in the the first sentence gives the result. 
\endpf
\smallskip \\

Notice the pseudohyperbolic metric appearing again in our calculations.
It is a fundamental artifact of geometric function theory (see [GAR]).

The following is an interesting and nontrivial fact about the 
Carath\'{e}odory and Kobayashi metrics, one of the few
that is valid for a large class of domains.

\begin{theorem} \sl
Let $\Omega$ be any planar domain that is not simply connected.  Then the Kobayashi metric and
the Carath\'{e}odory metric are unequal on $\Omega$.
\end{theorem}

\noindent {\bf Remark:}  We shall prove this result {\it not} by actually calculating the
metrics, but rather by an indirect argument.  It is a pleasing application of
geometric analysis.
\endpf
\smallskip \\

We shall first need a lemma.
\smallskip   \\

\begin{lemma}  \sl
Let $\Omega$ be as in the theorem and $D$ the unit disc as usual.  Then there do {\it not}
exist holomorphic functions $\varphi: D \ra \Omega$ and $\psi: \Omega \ra D$ such
that $\psi \circ \varphi(z) \equiv z$.
\end{lemma}
{\bf Proof:}  We thank Paul Gauthier for this elegant proof.  

Obviously the function $\varphi$ must be one-to-one, otherwise
the composition could not be one-to-one.  We claim that
$\varphi$ is onto.  If that is the case, then $\varphi$ is a conformal
equivalence and hence, in particular, a homeomorphism.  But the disc
and $\Omega$ cannot be homeomorphic (their first homotopy groups are 
different, for example).

To prove the claim, suppose not.  Then the image of $\varphi$ is a proper open
subset of $\Omega$.  Call the image $U = \varphi(D)$.  Then $U$ has a boundary
point $p \in \Omega$.  See Figure 3, which illustrates the idea when $\Omega$ is
an annulus.  Since $\varphi$ is continuous, it takes
compact sets to compact sets.  It follows that $\varphi^{-1}$ has the
property that if $U \ni z_j \ra p$ then $w_j \equiv \varphi^{-1}(z_j)$ tends to
some boundary point $b$ of $D$.  But then $w_j = \psi \circ \varphi(w_j) = \psi(z_j)$
for each $j$.  Letting $j \ra \infty$, we see that 
$$
b = \lim_{j \ra \infty} w_j = \lim_{j \ra \infty} \psi(z_j) = \psi(p) \equiv q \in D \, .
$$
Thus we have the boundary point $b$ of $D$ equaling a interior point $q \in D$.  This
is impossible.
\endpf
\smallskip \\

    \begin{figure}
    \centering
      \includegraphics[height=1.75in, width=5.3in]{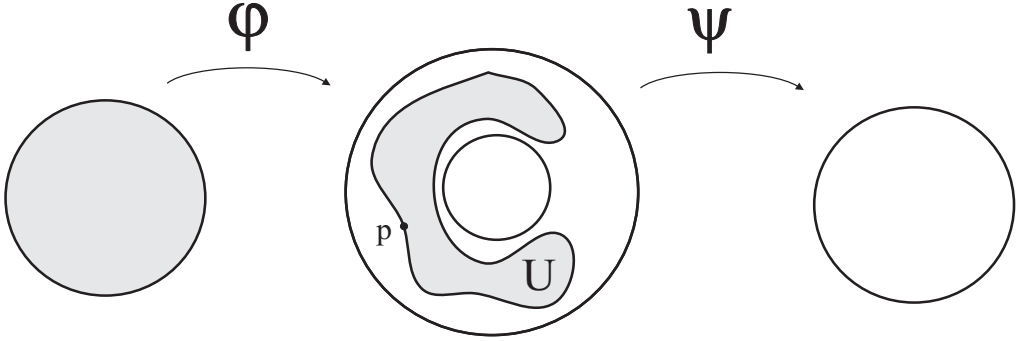}
      \caption{The identity on the disc does not factor through $\Omega$ (the case of the annulus).}
    \end{figure}

\noindent {\bf Proof of the Theorem:}  Fix a point $p \in \Omega$.  We claim that there is a constant $0 < c < 1$ such
that if $\varphi: D \ra \Omega$ is holomorphic with $\varphi(0) = p$ and $\psi: \Omega \ra D$ is
holomorphic with $\psi(p) = 0$ then 
$$
\left | (\psi \circ \varphi)'(0) \right | \leq c \, .	\eqno (\star)
$$

Suppose not.  Then, for each integer $j > 0$, there are holomorphic $\varphi_j: D \ra \Omega$
and $\psi_j: \Omega \ra D$, with $\varphi_j(0) = p$ and $\psi_j(p) = 0$, such that
$$
\left | (\psi_j \circ \varphi_j )'(0) \right | > 1 - \frac{1}{j} \, .
$$
Applying Montel's theorem, we may extract subsequences $\varphi_{j_k} \ra \varphi_0$
(uniformly on compact sets) and $\psi_{j_k} \ra \psi_0$ (uniformly on
compact sets).  Of course it will be the case that $\varphi_0: D \ra \Omega$, $\psi_0: \Omega \ra D$,
$\varphi_0$ and $\psi_0$ are holomorphic, and $\varphi_0(0) = p$, $\psi_0(p) = 0$.
And, what is most important,
$$
\left | (\psi_0 \circ \varphi_0)'(0) \right | = 1 \, .
$$
By Schwarz's lemma, we may conclude that $\psi_0 \circ \varphi_0$ is a rotation.
Postcomposing $\psi_0$ with the inverse of that rotation, we end up 
with a map from the disc to $\Omega$ and another map from $\Omega$
to the disc so that their composition from the disc to the disc is the identity.
The lemma tells us that this is impossible.

Now inequality $(\star)$ tells us that
$$
|\psi'(p)| \leq c \cdot \frac{1}{|\varphi'(0)|} \, .
$$
Taking the infimum of the righthand side over $\varphi$ and the supremum
of the lefthand side over $\psi$ as usual yields that
$$
F_K^\Omega(p, 1) \leq c \cdot F_C^\Omega(p, 1) \, .
$$
That is the desired result.
\endpf
\smallskip \\

We conclude this section by recording an interesting fact about our two invariant metrics.
This will prove useful in the applications presented in the next section.  The result
has been somewhat anticipated in our discussion of the boundary behavior of the Kobayashi metric
on the unit disc.

\begin{proposition} \sl
Let $\Omega \subseteq \CC$ be a bounded domain with $C^2$ boundary (i.e., the boundary is locally
the graph of a $C^2$ function).  Then there are constants $c, C > 0$ such that, with
$\delta(z)$ denoting the Euclidean distance to the boundary of $\Omega$, 
$$
\frac{c|\xi|}{\delta(z)} \leq F_K^\Omega (z, \xi) \leq \frac{C|\xi|}{\delta(z)} \, .
$$
A similar set of estimates holds for the infinitesimal \cara\ metric.  
\end{proposition}

\noindent {\bf Remark:}	  A glance at the proof of the proposition shows that the upper bound is true for {\it any}
domain that is a proper subset of $\CC$.  For the disc $D(z, \delta(z))$ certainly lies
in $\Omega$, and then elementary comparison (Corollary 5) gives the result.

In the language of Gehring and Palka [GEP], Proposition 17 shows that the Kobayashi (or Carath\'{e}odory)
metric is comparable to the quasihyperbolic metric.
\endpf
\smallskip \\				  

\noindent {\bf Proof of the Proposition:}  It follows from the idea of the osculating circle (see [BLK]) in multivariable
calculus that there is are numbers $r, R > 0$ such that each point $p \in \partial \Omega$ 
has an {\it interior} osculating circle $C(p',r)$ at $p$ and an {\it exterior} osculating
circle $C(p'', R)$ at $p$.  See Figure 4.

Now if $z \in \Omega$ and $z$ is sufficiently near the boundary then, by the tubular
neighborhood theorem (see [HIR]), there is a unique nearest point $\pi(z) \in \partial \Omega$.
Consider the osculating circle $C(\pi(z)', r)$ at that point and its corresponding
disc $D(\pi(z)', r)$.  Then, by Corollary 5,
$$
F_K^\Omega(z, \xi) \leq F_K^{D(\pi(z)', r)} (z, \xi) \approx \frac{c|\xi|}{\delta(z)} \, .
$$
That is one half of what we wish to prove.  [Notice that we need only consider $z$ near the
boundary since the estimates are trivial in the interior.]

    \begin{figure}
    \centering
      \includegraphics[height=2.65in, width=2.75in]{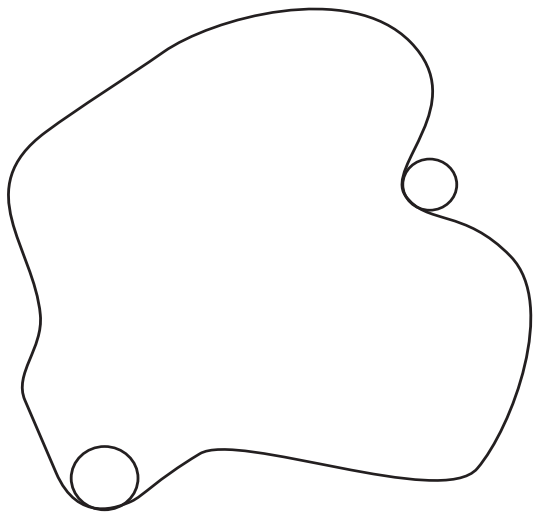}
      \caption{Interior and exterior osculating circles.}
    \end{figure}

For the other half, we again consider $z \in \Omega$, $z$ near the boundary.  Again let 
$\pi(z)$ be the nearest point in the boundary.  Let $D(\pi(z)'', R)$ be the exterior
osculating disc at $\pi(z)$.  Let $D(\pi(z)'', \widetilde{R})$ be a large disc
centered at $\pi(z)''$ that contains the domain $\Omega$.  Now consider the
region $U \equiv D(\pi(z)'', \widetilde{R}) \setminus \overline{D(\pi(z)'',R)}$.
Then certainly $U \supseteq \Omega$.  Hence, by Corollary 5,
$$
F_K^\Omega(z, \xi) \geq F_K^U(z, \xi) \, .
$$
But we may use a simple ``reflection'' map $\zeta \mapsto R^2/(\zeta - \pi(z)'')$ to
compare the Kobayashi metric on $U$ with the Kobayashi metric on a disc and then see that 
$$
F_K^U (z, \xi) \approx \frac{c |\xi|}{\delta(z)} \, .
$$
Putting together the last two displayed lines yields
$$
F_K^\Omega (z, \xi) \geq \frac{c|\xi|}{\delta(z)} \, .   \eqno \BoxOpTwo
$$

It is an immediate corollary of Proposition 17 that the {\it Euclidean} diameter\footnote{The Euclidean
diameter of a set is the supremum of Euclidean distances of pairs of points in the set.}
of a Carath\'{e}odory or Kobayashi metric ball $B(p,r)$ (for $r > 0$ fixed) tends to 0 as $p$ tends
to the boundary of a $C^2$ bounded domain.  We leave the details of this assertion
for the interested reader.

Certainly one important upshot of Proposition 17 is that, on a domain with $C^2$ boundary,
the Kobayashi metric is complete (a similar assertion holds for the Carath\'{e}odory metric).
That is so because we know that the metric blows up like the reciprocal of the distance to the
boundary.  Thus we can see, with some tedious but straightforward calculations (just
as we did on the unit disc), that the length of any curve tending to the boundary
is infinite.

\section{Some Applications}

We now show how metric geometry can in fact inform our study of function theory.\footnote{As noted
earlier, the roles of the Carath\'{e}odory and the Kobayashi metrics are essentially
interchangeable in these examples.}
The first result is due to Farkas and Ritt, but the proof is due
to Earle and Hamilton [EAH].  It concerns fixed points for holomorphic functions.
It is pleasing because it uses not only one of our invariant metrics, but it
also uses a fixed point-theorem from functional analysis.

\begin{theorem}[Farkas, Ritt]  Let  
$f: D \rightarrow D$  be holomorphic and assume that the image  
$M = \{ f(z): z \in  D \}$  of  $f$ has compact closure in  $D.$
Then
there is a unique point  $P \in  D$  such that  $f(P) = P.$
We call  $P$  a {\it fixed point} for  $f.$ 
\end{theorem}

\noindent {\bf Remark:}  This result is actually amenable to a number
of different proofs.  
If we take the image of $f$ to lie in a disc $\overline{D}(0, r)$ for some
$0 < r < 1$ then we may think of $f$ as mapping $\overline{D}(0,r)$ to
$\overline{D}(0,r)$ continuously.  Thus the Brouwer fixed-point theorem
applies and we find a fixed point (although this argument does not address
the uniqueness question).

It is well known that {\it any} proof of a fixed point theorem
will involve argument principle considerations (that is what
homotopy does for us in Brouwer's original proof; see also the
proof in [GAG]). Thus it stands to reason that Rouch\'{e}'s
theorem can be used to give the present result. That approach
also does not address the uniqueness question. Our purpose
here is to illustrate the utility of metric geometry, and also
to derive the stronger uniqueness result for the fixed point.
\endpf
\smallskip \\

\noindent {\bf Proof:}  By hypothesis, there is an  $\epsilon > 0$  such that if 
$m  \in   M$  and   $|  z |  \geq   1$  then  $|  m - z |   > 2
\epsilon$.  See Figure 5.  Fix  $z_0  \in D$  and define 
$$g(z) = f(z) + \epsilon (f(z) - f(z_0)).$$
\noindent Then  $g$  is holomorphic and  $g$  still maps  $D$  into
$D.$ Also
$$
g'(z_0) = (1 + \epsilon) f'(z_0).
$$
Certainly $g$ is distance-decreasing
in the Kobayashi metric. Therefore		       
$$
 F_K^D (g(z_0),g'(z_0) \cdot \tau) \leq  F_K^D (z_0, \tau)
$$
for any tangent vector $\tau$.  Now if  $\gamma: [a,b] \rightarrow D$  
is any continuously differentiable curve, and if we take $t \in [a,b]$, 
$z_0 = \gamma(t)$, and $\tau = \gamma'(t)$, then we may conclude that
$$
F_K^D (g(\gamma(t)), g'(\gamma(t)) \cdot \gamma'(t)) \leq F_K^D(\gamma(t), \gamma'(t)) \, .
$$
Writing this out gives
$$
(1 + \epsilon) F_K^D(f(\gamma(t)), f'(\gamma(t)) \cdot \gamma'(t)) \leq F_K^D (\gamma(t), \gamma'(t)) \, .
$$

Integrating both sides from $a$ to $b$, we conclude that
$$
L_K^D  (f \circ \gamma ) \leq  (1 + \epsilon )^{-1} L_K^D (\gamma ) \, .
$$

    \begin{figure}
    \centering
      \includegraphics[height=1.95in, width=4.3in]{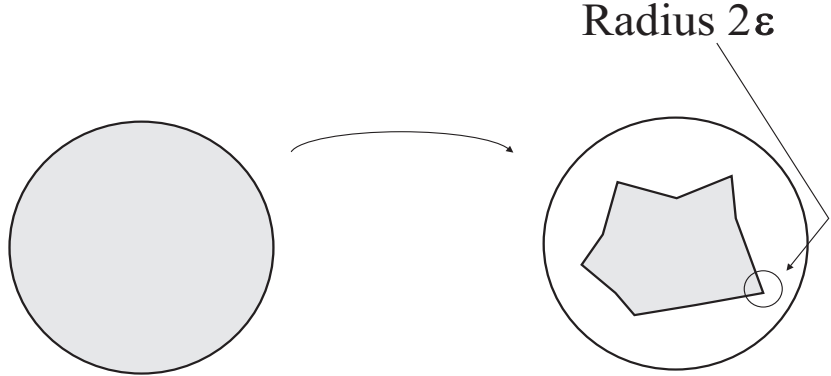}
      \caption{Relative compactness of the image of the mappings.}
    \end{figure}

\noindent If  $P$  and  $Q$  are elements of  $D$ then we see that
$$
K (f(P),f(Q))  \leq  (1 + \epsilon )^{-1} K (P,Q) \ .
$$

We conclude that  $f$  is a contraction in the Kobayashi metric.
Recall
that in Section 6 we proved that the disc  $D$  
is a complete metric space when equipped with the Kobayashi
metric.  By the contraction mapping fixed-point theorem (see [LS]), 
$f$  has a unique fixed point. 
\endpf 
\smallskip \\

Now we shift gears and look at the boundary behavior of holomorphic 
functions.  Complete background may be found in [KRA2, Ch.\ 8] or [KRA3].  We begin
by reviewing some terminology.  Let $f$ be a function (not necessarily
holomorphic) on the unit disc $D$.  Let $p = e^{i\theta}$ be a point
in the boundary of the disc.   We say that $f$ has {\it radial boundary
limit} $\ell$ at $p$ if
$$
\lim_{r \ra 1^-} f(r p) = \ell  \, .
$$
As a counterpoint, let us now consider a broader notion of limit.
For $p = e^{i\theta} \in \partial D$ and $\alpha > 1$, let us 
define
$$
\Gamma_\alpha (p) = \{z \in D: |z - p| < \alpha (1 - |z|) \} \, .
$$
See Figure 6.	In fact an analogous definition works just as well on any
domain with $C^2$ boundary (with $(1 - |z|)$ replaced by $\delta(z)$, the
distance of $z$ to the boundary).

    \begin{figure}
    \centering
      \includegraphics[height=1.8in, width=4.3in]{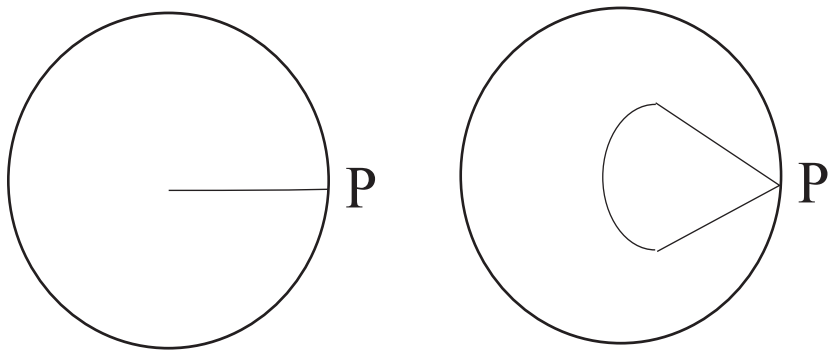}
      \caption{Radial convergence and nontangential convergence.}
    \end{figure}

We say that $f$ on $D$ has {\it nontangential} limit $\ell$ at $p$ if 
$$
\lim_{\Gamma_\alpha(p) \ni z \ra p} f(z) = \ell 
$$
for each $\alpha > 1$.  Our purpose now is to compare and relate
radial convergence with nontangential convergence.  We shall work
on domains with $C^2$ boundary, which simply means that the boundary
is locally the graph of a twice continuously differentiable function.

\begin{theorem} \sl  Let  $\Omega_{1}, \Omega_{2}$  be bounded domains
with $C^{2}$  boundary and let
$$ 
f: \Omega_{1} \longrightarrow \Omega_{2} 
$$
\noindent be holomorphic.  If  $P \in \partial \Omega_{1},\ Q \in
\partial \Omega_{2} ,$ and  $f$  has radial limit  $Q$  at $P$ 
then $f$  has nontangential limit  $Q$  at $P$.	 
\end{theorem}

\noindent {\bf Remark:}  This is a version of a classical result
that is known as the {\it Lindel\"{o}f principle}.  The usual
proof of that result uses a normal families argument (see, for example,
[KRA2] or [KRA3]).  That argument is lurking in the background of the more geometric
argument that we present here.
\endpf
\smallskip \\

\noindent {\bf Proof:}  If  $z$  is an element of one of our domains  $\Omega_{j}$ 
and if  $s > 0$  then we let  $B(z,s)$  denote the metric
ball with center  $z$  and radius  $s$ in the Carath\'eodory
metric for $\Omega_{j}.$ For $P \in \partial \Omega_j$, we let
$\nu_P$ denote the unit outward normal at $P$ to the boundary $\partial \Omega_1$.
If $r_0 > 0$ and $\beta >0$ are fixed, we define
$$ 
{\cal M}_{\beta }(P) = \bigcup_{0<r< r_0}  B_{\Omega_1}(P - r \nu_{P},\beta) \, .
$$
Here $B$ denotes a metric ball.	 Observe that we use a subscript on $B$ 
to indicate in what domain the metric ball lives.

The estimate
$$ 
F^\Omega_C (z, \xi) \approx  {C |\xi| \over \hbox{dist}(z,\partial U)}\eqno(*)
$$
\noindent from Proposition 17 makes it a tedious but not difficult exercise to calculate
that the regions ${\cal M}_{\beta }$  are comparable to the regions 
$\Gamma_{\alpha }$ (in this last formula, and in what follows,
``dist'' means Euclidean distance).  
In point of fact, suppose that $z$ lies in
some $\Gamma_\alpha(p)$, some $p \in \partial \Omega_j$.
Let us denote by $\tau_p$ the {\it inward} normal {\it segment} ({\it not} the vector!) 
emanating from that boundary point $p$.  Using the estimate $(*)$, one can then
estimate that $\hbox{dist}_C (z, \tau_p) \leq C \cdot \alpha$.
For the converse estimate, assume that $z \not \in \Gamma_\alpha(p)$ 
and the same estimate shows that 
$\hbox{dist}_C (z, \tau_p) \geq C \cdot \alpha$.

Thus we see that 
$$ 
\lim_{\Gamma_\alpha
(P) \ni z \in P} f(z) = \ell , \qquad \forall \alpha > 1 
$$
\noindent iff
$$ 
\lim_{{\cal M}_\beta (P) \ni z \to P}  f(z) =  \ell , \qquad
\forall  \beta > 0 .\eqno(**) 
$$
\noindent Thus it is enough to prove  $(**).$

Since

$$ {\cal M}_{\beta } (P) = \bigcup_{0<r< r_0}  B_{\Omega_1}(P - r
\nu_{P},\beta ),
$$

\noindent the distance-decreasing property of  $f$  with
respect 
to the Carath\'eodory metric implies that
$$ 
f({\cal M}_{\beta } (P))  \subseteq \bigcup_{0<r< r_0}  B_{\Omega_2}(f(P - r
\nu_{P}),\beta) . 
$$

Pick  $\epsilon  > 0$.  By the radial limit hypothesis, 
there is a  $\delta  > 0$  such that if  $0 < t < \delta$ 
then
$$ 
| f(P - t\nu_{P}) - Q |   < \epsilon . 
$$
\noindent For such a $t,$ if $z \in B(P - t\nu_{P},\beta )$ then
$$ 
f(z)  \in   B_{\Omega_2}(f(P - t\nu_{P}),\beta ) . 
$$
\noindent But
$$ 
\hbox{dist}(f(P - t\nu_{P}),\partial \Omega_2) \leq 
\hbox{dist}(f(P - t\nu_{P}),Q) < \epsilon . 
$$
\noindent Therefore the estimate  $(*)$  implies that the
metric ball  $B(f(P - t\nu_{P}),\beta )$  has Euclidean diameter not
exceeding $C \cdot \epsilon.$  Here  $C$  depends on 
$\beta,$ but $\beta $  has been fixed once and for all.  Thus
$$ 
|f(z) - f(P - t\nu_{P}) |  < C \epsilon , \qquad 
\forall z \in B_{\Omega_1}(P - t\nu_P,\beta) . 
$$

We conclude that
\begin{eqnarray*}
|f(z) - Q| 
& \leq & |f(z) - f(P - t\nu_{P}) |  + |f(P - t\nu_{P}) - Q| \\
& \leq &  C \epsilon + \epsilon  = C' \epsilon \, . 
\end{eqnarray*}
\noindent This is the desired conclusion.
\endpf
\smallskip \\

Our final application concerns automorphism groups.  Some preliminary discussion
is in order.  If $\Omega$ is a planar domain then we consider the collection
of all conformal self-maps $\varphi: \Omega \ra \Omega$.  To be explicit, we demand
that $\varphi$ be holomorphic, one-to-one, and onto, and have a holomorphic
inverse.  This collection is a group when equipped with the binary operation
of composition of functions.  We call this the {\it automorphism group} of
$\Omega$, and we denote it by $\Aut(\Omega)$.

We endow the automorphism group with the topology of uniform
convergence on compact sets.  This is equivalent with the compact-open topology.
It is a fact that, with this topology, the automorphism group of a bounded domain
is a real Lie group (see [KOB2]).  We shall not need that information here.

One of the ways that we can understand a domain is by understanding its
automorphism group.  This may entail studying the group's algebraic properties, or
studying its topological properties, or perhaps by considering some
combination of the two.  The next result illustrates this symbiosis.

\begin{theorem} \sl  Let  $\Omega  \subseteq \CC$    be a
bounded domain with $C^{2}$  boundary.  If  $\hbox{\rm Aut}(\Omega)$  is
noncompact then  $\Omega$  is conformally equivalent to the unit
disc.
\end{theorem}

\noindent {\bf Remark:}  It is certainly known (see [MUR]) that a finitely
connected domain with connectivity at least 3 (i.e., at least two holes)
has only finitely many conformal self-maps.  This is a nontrivial result.  
Our approach gives another way to think about that classical result (which
was originally proved by Maurice Heins in the 1940s---see [HEI1], [HEI2]).
\endpf
\smallskip \\

\noindent We prove this theorem with a sequence of lemmas, 
each of which has intrinsic interest.  
\smallskip   \\

\begin{lemma} \sl Let  $\Omega  \subseteq \CC$    be bounded.
The group  $\hbox{\rm Aut}(\Omega)$  is compact if and only if, for
each 
$P \in \Omega$,  there is a compact  $K^{P} \subseteq  \Omega$  such that  
$\varphi (P)  \in  K^{P}$  for all  $\varphi  \in  \hbox{\rm
Aut}(\Omega).$
\end{lemma}
{\bf Proof:}  Assume that  $\hbox{\rm Aut}(\Omega)$  is compact.
Fix  $P  \in  \Omega.$  If there is no set  $K^{P}$  as claimed
then there exist  $\varphi_{j} \in \hbox{\rm Aut}(\Omega)$  such that 
$\varphi_{j}(P) \rightarrow w \in \partial \Omega,$  some  $w.$
But  $\Omega$  is bounded so that $\{ \varphi_{j} \}$  is a
normal family; thus there is a subsequence  $\varphi_{j_{k}}$  and a
holomorphic limit function $\varphi_{0}$  such that
$$ 
\varphi_{j_{k}} \longrightarrow \varphi_{0} 
$$
normally.

Notice that the image of each  $\varphi_{j}$  lies in  $\Omega$  
hence the image of  $\varphi_{0}$  lies in the  {\it closure}
$\overline{U}$  of $\Omega.$  If  $\varphi_{0}$  is nonconstant then
it satisfies the open mapping principle.  But
$$ 
\varphi_{0}(P) = \lim_{k \to \infty} \varphi_{j_{k}} (P) = w,
$$
\noindent hence the image of  $\varphi_{0}$  contains the
accumulation point  $w  \in  \partial \Omega,$ so it contains a neighborhood of 
$w.$ This is impossible because  $w$  is in the boundary of the
image of $\varphi_{0}.$  Therefore  $\varphi_{0}$  must be
constantly equal to  $w;$  thus  $\varphi_{0} \not \in  \hbox{\rm
Aut}(\Omega).$  The sequence  $\varphi_{j_{k}}$  therefore violates the
compactness of  $\hbox{\rm Aut}(\Omega).$ We conclude that  $K^{P}$  must
exist.

For the converse, fix  $P \in \Omega$  and let  $K^{P}$  
be the corresponding compact set in  $\Omega$  whose existence we
assume.  Let  $\{ \varphi_{j} \} \subseteq   \hbox{\rm
Aut}(\Omega)$ be any sequence.  Since $\Omega$  is bounded, there is a normally
converging subsequence  $\varphi_{j_{k}}$  with holomorphic
limit function $\varphi_{0}.$  As in the first half of the proof, 
if the image of  $\varphi_{0}$  contains any boundary point 
$w$ then  $\varphi_{0}$  must be constantly equal to  $w.$
But the image of  $P$  under  $\varphi_{0}$  must lie in 
$K^{P},$ so this possibility is ruled out.
We conclude that the image of  $\varphi_{0}$  lies in  $\Omega.$

Next notice that each $\varphi_{j_{k}}$ has an inverse $\psi_{j_{k}}.$
Passing to another subsequence, we may suppose that the $\psi_{j_{k}}$
converge to a limit function $\psi_{0}.$ For convenience, we denote this
final subsequence by $\psi_{m},$ corresponding to the automorphisms
$\varphi_{m}.$ Just as for $\varphi_{0},$ we can be sure that the image of
$\psi_{0}$ lies in $\Omega$. By Hurwitz's theorem, $\psi_0$ is
non-constant.

Now we have
$$ 
z \equiv \lim_{m \to \infty} \varphi_{m} \circ \psi_{m}(z)
= \varphi_{0} \circ \psi_{0}(z). 
$$
\noindent Since ${\bf i} (z) \equiv  z$  is onto, 
so is  $\varphi_{0}.$  Also, by the argument principle, the image of  $\psi_{0}$  is open, closed, and
nonempty.  Therefore $\psi_{0}$  is surjective.
Since  ${\bf i}  (z)$  is injective, it now follows that 
$\psi_{0}$  is injective.  Therefore $\psi_{0} \in  
\hbox{\rm Aut}(\Omega)$ and it follows that $\phi_0 \in \hbox{\rm Aut}(\Omega)$.
We conclude that
$$ 
\hbox{\rm Aut}(\Omega)  \ni \varphi_{m} \longrightarrow
\varphi_{0}
\in   \hbox{\rm
Aut}(\Omega), 
$$
\noindent and  $\hbox{\rm Aut}(\Omega)$  is compact.
\endpf 
\smallskip \\

\noindent {\bf Remark:}  It may be noted that, in the proof of the
converse direction of Lemma 21, only one compact set $K^P$ was needed.
\endpf
\smallskip \\

\begin{lemma} \sl Let  $\Omega  \subseteq  \CC$  be a bounded
domain with  $C^{2}$  boundary.  Suppose that  $P  \in  \Omega,\  \{
\varphi_{j} \}$  are holomorphic maps from  $\Omega$  to  $\Omega,$ and
$$ 
\varphi_{j}(P) \longrightarrow w  \in  \partial \Omega . 
$$
\noindent If  $K$  is compact in  $\Omega$  and  $V$  is a
neighborhood of  $w$  then there exists a positive number  $J$  such that
if $j \geq   J$  then
$$ 
\varphi_{j}(K)  \subseteq   V . 
$$
\noindent See Figure 7.
\end{lemma}

    \begin{figure}
    \centering
      \includegraphics[height=2.65in, width=2.75in]{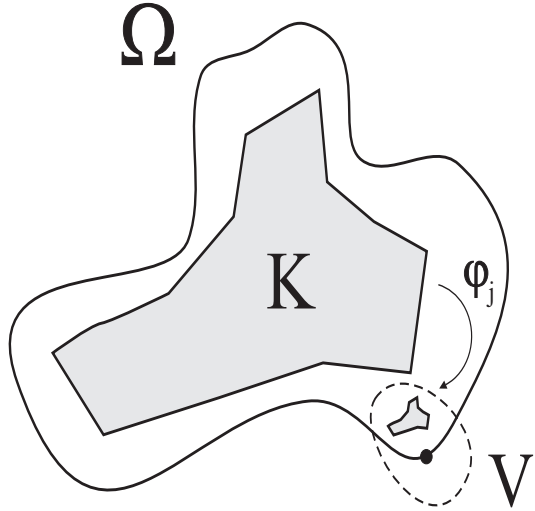}
      \caption{Noncompact group action on a compact set.}
    \end{figure}

\noindent {\bf Proof:}  Since  $\Omega,$ when equipped with the Carath\'eodory
metric, 
is a metric space and since  $K$  is compact, 
there is a positive number  $R$  such that the metric ball 
$B(P,R)$  contains  $K$.
Let  $Q_{j} = \varphi_{j}(P).$
Since each  $\varphi_{j}$  is distance-decreasing in the
Carath\'eodory metric, it follows that  $\varphi_{j}(B(P,R)) 
\subseteq  B(Q_{j},R).$  We claim that there is a positive 
$J$ 
such that whenever  $j  \geq   J$  then  $B(Q_{j},R) 
\subseteq  
V.$  Assuming the claim, we would then have
$$ 
\varphi_{j}(K)  \subseteq  \varphi_{j}(B(P,R)) \subseteq 
B(Q_{j},R) \subseteq  V, 
$$
as required.

To prove the claim, recall 
(because the Carath\'{e}odory metric
on $\Omega$ is complete) that the Euclidean radii of the metric balls 
$B(Q_{j},R)$  must tend to  $0.$  Choose  $\epsilon > 0$  such that the
Euclidean disc of center  $w$  and radius  $2\epsilon $  lies in  $V.$
We select  $J$  so large that when  $j > J$  then both the 
Euclidean distance of  $Q_{j}$  to  $w$  is less than 
$\epsilon$ and the Euclidean radius of  $B(Q_{j},R)$  is less than 
$\epsilon.$ The claim now follows from the triangle inequality.
\endpf
\smallskip \\

\noindent {\bf Proof of Theorem 20:} If  $\hbox{\rm Aut}(\Omega)$  is
not compact then, by Lemma 21, there is a sequence  $\varphi_{j}
\in \hbox{\rm Aut}(\Omega)$  and a  $P  \in  \Omega$  such that
$$ 
\varphi_{j}(P) \longrightarrow w  \in  \partial \Omega, 
$$
\noindent for some  $w  \in  \partial \Omega.$

Let
$$ 
\gamma: [0,1] \longrightarrow \Omega 
$$
\noindent be any continuous closed curve in  $\Omega.$
Since  $\partial \Omega$  is  $C^{2}$  there is a neighborhood 
$V$  of $w$  such that  $\Omega \cap V$  is simply connected 
(see Figure 8---the existence of the interior osculating
circle, or the tubular neighborhood, 
provided by the proof of Proposition 13 makes this assertion
clear).

    \begin{figure}
    \centering
      \includegraphics[height=2.65in, width=2.75in]{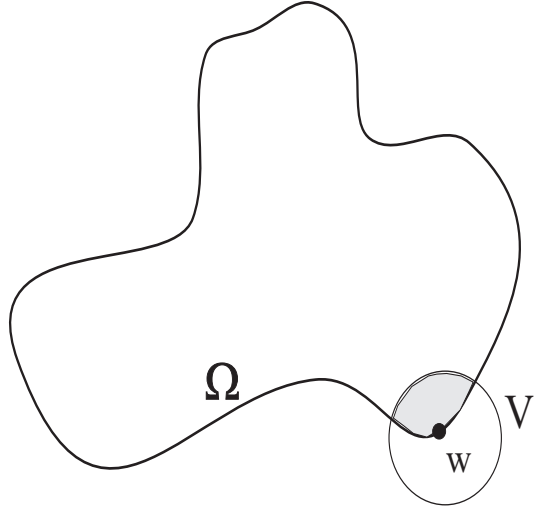}
      \caption{Simply connected boundary neighborhood.}
    \end{figure}

Let
$$ 
K = \{ \gamma (t): 0  \leq   t \leq   1 \} . 
$$
\noindent Then  $K$  is compact.
By Lemma 22, there is a  $J  \geq   0$  such that  $j  \geq  
J$ implies  $\varphi_{j}(K)  \subseteq   \Omega \cap  V.$
Thus  $\varphi_{j} \circ \gamma $  is a continuous, 
closed curve in  $\Omega \cap  V.$
The simple connectivity of  $\Omega \cap V$  implies that  
$\varphi_{j} \circ \gamma$  may be continuously deformed 
to the point  $\varphi_{j} \circ \gamma (0);$ that is,  
there is a homotopy
$$ 
\Psi : [0,1] \times [0,1] \longrightarrow \Omega \cap  V 
$$
\noindent such that
$$ 
\Psi (0,t) = \varphi_{j} \circ \gamma (t), \qquad \forall
t  \in [0,1] 
$$
\noindent and
$$ 
\Psi (1,t) = \varphi_{j} \circ \gamma (0), \qquad \forall t
\in [0,1] . 
$$
\noindent But then
$$ 
(\varphi_{j})^{-1} \circ  \Psi 
$$
\noindent is a homotopy of  the curve  $\gamma$  to the point 
$\gamma (0).$  It follows that $\Omega$ is simply connected.
By the Riemann mapping theorem, $\Omega$ is conformally equivalent to the disc.
\endpf
\smallskip \\
				
\section{Concluding Remarks}

Certainly the interaction of metric geometry with function theory has been
one of the seminal developments of twentieth century complex analysis.
There have many vectors in this activity:  {\bf (i)}  Poincar\'{e}, Bergman,
\cara, and Kobayashi (among several others) have provided us with
a family of extremely useful conformally invariant metrics, {\bf (ii)}  Lars
Ahlfors has shown that the Schwarz lemma may be understood in terms of curvature
of a suitable conformal metric (see [KRA1]),
and {\bf (iii)}  many of the phenomena of function theory have been given
very natural interpretations in terms of the geometry of K\"{a}hler manifolds.
Surely other authors would emend or modify this list.

The result of all these new ideas has been a subject enriched with new
results, and with new interpretations of old results.  Even the deep Picard
theorems may be given rather direct and quick proofs using metric geometry
(see [KRA1] for the details).  Each of the
applications presented in the present paper can actually be proved with
classical techniques.  But the metric geometry proofs are natural, enlightening,
and fun.

We hope that this excursion into the world of complex analysis and geometry
has provided the reader with adequate motivation to explore further.  The
result will be both edifying and rewarding.	     

\newpage

\noindent {\Large \sc References}
\smallskip \\

\begin{enumerate}

\item[{\bf [BER]}] S. Bergman, \"{U}ber Hermitesche unendliche Formen, die
zu einem Bereich geh\"{o}ren, nebst Anwendungen auf Fragen der Abbildung
durch Funktionen zweier komplexen Ver\"{a}nderlichen, {\it Math. Z.}
29(1929), 640--677.

\item[{\bf [BLK]}]  B. Blank and S. G. Krantz, {\it Calculus, Multivariable}, Key College Press,
Emeryville, CA, 2006.

\item[{\bf [CAR]}]  C. Carath\'{e}odory, \"{U}ber eine spezialle  Metrik, die in der Theorie der 
analytischen Funktionen auftritt, {\it Atti Pontifica Acad.\ Sc., Nuovi Lincei} 80(1927), 135--141.

\item[{\bf [DIF]}]  K. Diederich and J. E. Forn\ae ss, Comparison of the 
Bergman and Kobayashi metric, {\em Math. Ann.} 254(1980), 257-262.

\item[{\bf [EAH]}] C. Earle and R. Hamilton, A fixed point theorem for
holomorphic mappings, {\it Proc. Symp. Pure Math.}, Vol. XVI, 1968,
61--65.

\item[{\bf [FAK]}] H. Farkas and I. Kra, {\it Riemann Surfaces}, $2^{\rm
nd}$ ed., Springer-Verlag, New York, 1992.

\item[{\bf [FIS]}]  S. Fisher, {\it Function Theory on Planar
Domains}, John Wiley \& Sons, New York, 1983.

\item[{\bf [GAG]}]  T. W. Gamelin and R. E. Greene, {\it Introduction to Topology},
$2^{\rm nd}$ ed., Dover Books, Mineola, NY, 1999.

\item[{\bf [GAR]}]  J. Garnett, {\it Bounded Analytic Functions},
Academic Press, New York, 1981.

\item[{\bf [GEP]}]  F. Gehring and B. Palka, Quasiconformally 
homogeneous domains, {\it J. Analyse Math.} 30(1976), 172--199.

\item[{\bf [GOL]}] Goluzin, {\it Geometric Theory of Functions of a Complex
Variable}, American Mathematical Society, Providence, 1969.

\item[{\bf [GRK]}]  R. E. Greene and S. G. Krantz, {\it Function Theory of
One Complex Variable}, $3^{\rm rd}$ ed., American Mathematical Society,
Providence, RI, 2006.

\item[{\bf [HEI1]}]  M. Heins, A note on a theorem of Rad\'{o} concerning
the $(1,m)$ conformal maps of a multiply-connected region into
itself, {\it Bull.\ Am.\ Math.\ Soc.} 47(1941), 128--130.

\item[{\bf [HEI2]}]  M. Heins, On the number of 1--1 directly conformal
maps which a multiply-connected plane region of finite connectivity
$p$ ($>2$) admits onto itself, {\it Bull.\ Am.\ Math.\ Soc.} 52(1946),
454--457.

\item[{\bf [HIR]}]  M. Hirsch, {\it Differential Topology}, Springer-Verlag,
New York, 1976.

\item[{\bf [HUS]}]  D. Husemoller, {\it Fibre Bundles},
$2^{\rm nd}$ ed., Springer-Verlag, New York, 1975.

\item[{\bf [JAP]}]  M. Jarnicki and P. Pflug, {\it Invariant Distances and Metrics in
Complex Analysis}, de Gruyter, Berlin and New York, 1993.

\item[{\bf [KOB1]}]  S. Kobayashi, Invariant distances on complex manifolds and holomorphic
mappings, {\it J. Math.\ Soc. Japan} 19(1967), 460--480.

\item[{\bf [KOB2]}]  S. Kobayashi, {\it Hyperbolic Manifolds and Holomorphic
Mappings}, Dekker, New York, 1970.	

\item[{\bf [KOB3]}]  S. Kobayashi, {\it Hyperbolic Complex Spaces}, Springer-Verlag,
New York, 1998.

\item[{\bf [KRA1]}]  S. G. Krantz, {\it Complex Analysis:  The Geometric Viewpoint},
$2^{\rm nd}$ ed., Mathematical Association of America, Washington, D.C., 2004.
	
\item[{\bf [KRA2]}]  S. G. Krantz, {\it Function Theory of Several Complex Variables},
$2^{\rm nd}$ ed., American Mathematical Society, Providence, RI, 2001.

\item[{\bf [KRA3]}] S. G. Krantz, {\it Cornerstones of Geometric Function
Theory: Explorations in Complex Analysis} Birkh\"{a}user, Boston, 2006.

\item[{\bf [LS]}]  L. Loomis and S. Sternberg, {\it Advanced Calculus},
Jones \& Bartlett, Boston, MA, 1990. 

\item[{\bf [MUR]}]  C. Mueller and W. Rudin, Proper holomorphic self-maps of plane regions,
{\it Complex Variables} 17(1991), 113--121.

\item[{\bf [ROY]}] H. Royden, Remarks on the Kobayashi Metric, {\it Several
Complex Variables II}, Maryland 1970, Springer, Berlin, 1971, 125-137.

\item[{\bf [SIM]}] R. R. Simha, The Carath\'{e}odory metric of the annulus,
{\it Proc.\ AMS} 50(1975), 162--166.

\item[{\bf [SPA]}] E. Spanier, {\it Algebraic Topology}, Springer-Verlag,
New York, 1981.

\end{enumerate}

\end{document}